\documentclass{scrartcl}

\usepackage{hyperref}
\usepackage[english]{babel}
\usepackage[thms]{packages}
\usepackage{todonotes}
\usepackage{caption}
\usepackage{subcaption}

\usepackage{dsfont}

\usepackage[style=trad-alpha, backend=biber, eprint=false, isbn=false]{biblatex}

\title{The Skorohod Topologies}
\subtitle{What They Are and Why We Should Care}
\author{Julian Kern}

\begin{document}
\maketitle

\begin{abstract}
This paper presents a gentle and informal introduction to the Skorokhod topologies. Focus is on motivating examples and concepts.
\end{abstract}

\tableofcontents
\newpage

\section{Introduction}
\label{intro}
As soon as we stray from the theory of continuous stochastic processes, we are in need of a suitable space of discontinuous functions and a topology on it. 
Skorokhod proposed in \cite{Sko} the topology used predominantly today and which has since inherited his name.
When I started to work with discontinuous stochastic processes and needed to understand the Skorokhod space, I struggled to find textbooks or lecture notes providing an easy start into the topic. The general tenor is that ``constructing [the] Skorokhod topology and deriving tightness criteria are rather tedious" (see \cite[Chapter VI]{JS03}). That gave me the impression that the Skorokhod topology is a very technical tool which has no real motivation. 

After working with it for some years, I believe that there are simple and intuitive ideas \emph{underlying} this construction which might facilitate the understanding. Unfortunately, these are not the main focus in most textbooks as the proofs are already long enough as is. For the same reason, very few textbooks explore all four Skorokhod topologies and focus only on the \emph{main} one, also known as $J_1$-topology. Here, I want to take the time to expose the ideas underlying the four Skorokhod topologies. 

This paper is not meant to give a complete overview on the Skorokhod topologies, nor will I include any proofs. Instead, I will concentrate on pictures, examples and heuristics in the hope of building intuition. Nevertheless, I try to give a reference to every major fact I mention. 
Note that I will not discuss non-Skorokhod topologies such as the $S$-topology defined in \cite{J97}.

Most of what I present is taken from the two books \cite{Billingsley,SPL} which I highly recommend for delving into the subject: the former focuses on the prevalent topology on real-valued processes, the latter takes a more general approach. Some general results are taken from \cite{EK86}. For examples on how to apply these general results to measure-valued processes, I recommend the first chapter of \cite{Etheridge}.

The rest of the paper has a very simple structure: first, I try to motivate why we need a new topology and what we should expect it to look like. Next, I derive the two main $J_1$- and $M_1$-topologies for real-valued processes on a finite time interval. Eventually, I conclude by presenting generalisations of these topologies.

\paragraph*{Acknowledgement}

I would like to thank Andrey Pilipenko for bringing the collection \cite{SW} of selected works of Skorokhod to my attention. This research has been partially funded by Deutsche Forschungsgemeinschaft (DFG) through grant CRC 1114 Scaling Cascades in Complex Systems, Project
Number 235221301, Project C02 Interface dynamics: Bridging stochastic and hydrodynamic descriptions.

\section{Motivation}

In this first section, I illustrate why we need to define a new topology and what properties we would want it to have.

\subsection{Convergence of Continuous Processes}

Let us start with the most famous example of convergence of stochastic processes: Donsker's Theorem, a.k.a.~the functional central limit theorem. Consider a simple random walk $(S_n)_{n\geq 0}$ on $\mathbb{Z}$ starting at $S_0 = 0$, and define the rescaled and interpolated continuous process
\[
Y^N_t := \dfrac{1}{\sqrt{N}}\Big{(} S_{\lfloor Nt\rfloor} + \left(Nt - \lfloor Nt\rfloor\right)\left(S_{\lfloor Nt\rfloor +1} - S_{\lfloor Nt\rfloor}\right)\Big{)},\quad t\in [0,1].
\]
We consider $(Y^N_\cdot)_{N\geq 1}$ as a random sequence in $\mathcal{C}([0,1])$ endowed with the Borel $\sigma$-algebra of the topology of uniform convergence. Donsker's Theorem states that the sequence $(Y^N_\cdot)_{N\geq 1}$ converges in distribution to a standard Brownian motion $B$ on the time interval $[0,1]$. To prove this, one uses that $\mathcal{C}([0,1])$ endowed with the topology of uniform convergence is a Polish space\footnote{A \emph{Polish space} is a topological space that is separable and completely metrizable.} so that we may apply Prokhorov's Theorem.
\begin{theo}[Prokhorov; see {\emph{e.g.}~\cite[Section 5]{Billingsley} or \cite[Theorems 8.6.2 and 8.9.4]{B06} for more complete statements}]\label{theo:prokhorov}
Let $E$ be a Polish space with its Borel $\sigma$-algebra. Let $\mathcal{P}(E)$ denote the set of probability measures on this measurable space endowed with the topology of weak convergence of measures. Then the following holds true:
\begin{enumerate}
	\item the space $\mathcal{P}(E)$ is again Polish,
	\item a set of probability measures $K\subseteq \mathcal{P}(E)$ is relatively compact if and only if it is tight\footnote{A family $(P_\alpha)_\alpha$ of probability measures is \emph{tight} if they vanish uniformly outside compact sets: for every $\epsilon > 0$ there is a compact set $K_\epsilon$ such that $\sup_\alpha P_\alpha(K_\epsilon^c) <\epsilon$.}.
\end{enumerate}
\end{theo} 

\begin{corol}\label{corol:convergence_law_equiv}
A sequence of processes $(X^N)_{N\geq 1}\subset \mathcal{C}([0,1])$ converges in law to some $X\in \mathcal{C}([0,1])$ if and only if the sequence $(X^N)_{N\geq 1}$ is tight and the finite-dimensional distributions of $(X^N)_{N\geq 1}$ converge to those of $X$, \emph{i.e.}~for all $0\leq t_1 < \dots < t_k \leq 1$, one has
\[
\big{(}X^N_{t_1},\dots, X^N_{t_k}\big{)}\overset{(d)}{\underset{N\to +\infty}{\longrightarrow}} \big{(} X_{t_1},\dots, X_{t_k}\big{)}
\]
in $\mathbb{R}^k$.
\end{corol}

That means that one only needs to check those two things to prove Donsker's Theorem. More precisely, we really need to worry only about tightness as the convergence of the finite-dimensional distributions is an application of the usual Central Limit Theorem. Since tightness is a statement about compact sets of the underlying space, we need a good characterisation of compact sets in $\mathcal{C}([0,1])$. For this, we use another powerful theorem:
\begin{theo}[Arzelà--Ascoli; see {\emph{e.g.}~\cite{GV61}}]\label{theo:arzela-ascoli}
A set $K\subseteq \mathcal{C}([0,1])$ of continuous functions is relatively compact if and only if
\begin{enumerate}[label=\roman*)]
	\item the set $K$ is uniformly bounded at $0$ in the sense that
	\[
	\sup_{f\in K} \vert f(0)\vert < +\infty,
	\]
	\item the set $K$ is uniformly equicontinuous, \emph{i.e.}~for every $\epsilon > 0$ there exists some $\delta > 0$ such that
	\[
	\vert f(t) - f(s)\vert < \epsilon
	\]
	for all $f\in K$ whenever $\vert t-s\vert < \delta$. 
\end{enumerate}
\end{theo}
In terms of the modulus of continuity
	\[
	\omega_\delta(f) := \sup_{\vert t-s\vert < \delta} \vert f(t) - f(s)\vert,
	\]
the second condition may be rewritten as
	\[
	\lim_{\delta\downarrow 0}\sup_{f\in K} \omega_\delta(f) = 0.
	\]
\begin{corol}
A sequence of processes $(X^N)_{N\geq 1}\subset \mathcal{C}([0,1])$ is tight if and only if 
\begin{enumerate}[label=\roman*)]
	\item the sequence is bounded at $0$ in the sense that for all $\epsilon > 0$, one has
	\[
	\lim_{C\to +\infty}\sup_{N\geq 1} \mathbb{P}(\vert X^N_0\vert > C) < \epsilon,
	\]
	\item for every $\epsilon > 0$, one has
	\[
	\lim_{\delta\downarrow 0}\sup_{N\geq 1} \mathbb{P}\left[ \omega_{X^N}(\delta) > \epsilon\right] = 0.
	\]
\end{enumerate}
\end{corol}

With all this machinery at our disposal, it becomes straightforward to prove Donsker's Theorem. One only needs to do two things:
\begin{enumerate}
	\item use Prokhorov's Theorem to prove relative compactness via tightness;
	\item identify all limit points through their finite-dimensional distributions.
\end{enumerate}
Don't get me wrong: each point in itself might be difficult to prove. But at least we have a \emph{strategy} on how to approach the problem. It turns out that this strategy is very general: only the first point depends on the particular topology which we put on the space of functions. For it to work, we first and foremost need Prokhorov's Theorem. That means that we want a Polish topology. But there is a second ingredient on which we relied heavily. To prove tightness, we need a \emph{good way to characterise the compact sets} of the topology we work in. In the case of the uniform topology on $\mathcal{C}([0,1])$, this is taken care of by the Arzelà--Ascoli Theorem.

\subsection{Convergence of Discontinuous Processes}\label{ssec:disc_proc}

We will now try to apply our insights from the previous section to the convergence of discontinuous processes. More precisely, we first need to identify what space of functions we are interested in and what topology we can endow it with.

A first idea could be to go to the next bigger space we are familiar with and which extends the topology of uniform convergence: the space of bounded measurable functions $\mathcal{B}([0,1])$ with the topology of uniform convergence. The way I am presenting this, it becomes clear that this is not a good choice: even though the space is complete, it is not separable and therefore not Polish. And when we want to do probability theory, that is not a good sign; particularly with Prokhorov's Theorem in mind. Despite having quite a nice characterisation of the compact sets of $\mathcal{B}([0,1])$ similar to the Arzelà--Ascoli Theorem \ref{theo:arzela-ascoli}, it is not the right space to work in.

Now that we ruled out the obvious choice, we need to decide on how to proceed. The first step is to choose the right space of functions. In other words: what type of functions are relevant to us? Note that we are mostly interested in martingales and Markov processes, often characterised through their generator or their martingale problem. For these, we have a very nice regularity result:
\begin{theo}[see {\emph{e.g.}~\cite{L09}}]\label{theo:cadlag}
A sub- or supermartingale $(M_t)_{t\geq 0}$ with respect to a right continuous filtration has a \emph{càdlàg}\footnote{The acronym \emph{càdlàg} comes from the French ``continue à droite, limite à gauche" which translates to ``continuous from the right with left limits" (at any point $t$).} modification whenever $t\mapsto \mathbb{E}[M_t]$ is right continuous. In particular, a martingale w.r.t.~a right continuous filtration always has a \emph{càdlàg} modification.
\end{theo}

That indicates that the space of \emph{càdlàg} functions seems to be the right choice. We will denote this so-called \emph{Skorokhod space}\footnote{I have the impression that there are differences in the nomenclature. Sometimes, the name \emph{Skorokhod space} is only used for the space endowed with the usual \emph{Skorokhod topology} which we will construct later on. However, in different contexts other \emph{Skorokhod topologies} are useful, so that it is necessary to give different names to the \emph{topological} spaces.} by $\mathbb{D}_{[0,1]}(\mathbb{R})$. In general, the Skorokhod space of \emph{càdlàg} functions on $[0,T]$ (resp.~$[0,+\infty)$) with values in a (hopefully Polish) space $E$ will be denoted by $\mathbb{D}_{[0,T]}(E)$ (resp.~$\mathbb{D}_{[0,+\infty)}(E)$). Note that we may view $\mathcal{C}([0,1])$ as a subspace of $\mathbb{D}_{[0,1]}(\mathbb{R})$.

Now that we have identified the ``right" space of functions, we need to identify the ``right" topology. It turns out that there is not \emph{one} good topology. So instead, we will identify the right properties a good topology should have. The most important part is the applicability of Prokhorov's Theorem. In other words, we want
\begin{equation}\label{eq:polish}
	\text{The topology is Polish.}\tag{Polishness Property}
\end{equation}
to hold. In most applications, it is enough to weaken this condition to 
\begin{equation}\label{eq:prokhorov}
\begin{array}{cc}
	\text{The space is separable; and if a family of measures on it is tight,}\\\text{then it is relatively compact w.r.t.~the topology of weak convergence.}\tag{SepProkhorov Property}
\end{array}
\end{equation}
which amounts to the ``important" part of Prokhorov's Theorem. However, it is usually preferable to have a Polish space. The second important ingredient in our strategy was the Arzelà--Ascoli Theorem that characterises the compact sets for the topology of uniform convergence. Hence, we would want
\begin{equation}\label{eq:aa}
	\text{An Arzelà--Ascoli type theorem describing compact sets exists.}\tag{ArzAsc Property}
\end{equation}
to hold.

If both conditions \eqref{eq:polish} and \eqref{eq:aa} are satisfied, we have a ``good" topology. Nevertheless, there are other properties that one could wish for. For example, it would be great if the new topology extends the topology of uniform convergence on $\mathcal{C}([0,1])$. More formally, this means that
\begin{equation}\label{eq:extension}
	\text{The trace\footnotemark\ topology on $\mathcal{C}([0,1])$ is the topology of uniform convergence.}\tag{Extension Property}
\end{equation}\footnotetext{{The \emph{trace} topology, also known as \emph{subspace} topology, is obtained by restricting all open sets to the subspace.}}
Even though this property seems very natural, there is an important argument against it: the space of continuous functions with the topology of uniform convergence is \emph{complete}. That means that if the new topology extends it, it is impossible for a sequence of continuous functions to converge to a discontinuous function. In other words, a topology extending the topology of uniform convergence may be too strong for some applications. 

There is a last ``bonus" property which would be nice to have. Imagine a sequence of functions converging to some \emph{continuous} function. Since the limit lies in the subspace, where the topology is ``stronger", it would be great if this would automatically strengthen the mode of convergence to uniform convergence, \emph{i.e.~}
\begin{equation}\label{eq:BONUS}
	\text{If $X_n \rightarrow X$ and $X$ is continuous, then $X_n \rightarrow X$ uniformly.}\tag{Bonus Property}
\end{equation}
If the new topology satisfies \eqref{eq:BONUS}, we would not need to worry about the topology of uniform convergence any more at all. Whenever the limit is continuous, we get the uniform convergence for free!\\

Equipped with these constraints, we will start constructing Skorokhod topologies!

\section{Skorokhod Topologies on $\mathbb{D}_{[0,1]}(\mathbb{R})$}

In this main section, we construct the topologies on the Skorokhod space: first, we will try to understand how we might want to tweak the uniform topology; then, we will see the actual definitions of these topologies. A little word of caution: the rigorous proofs of the facts that I will state are very technical. So instead, I will use the old magician's trick and refer to the books \cite{Billingsley,SPL} that present those proofs nicely.

\subsection{What Exactly Goes Wrong with the Uniform Topology?}

In Section \ref{ssec:disc_proc}, I pointed out that one major problem of the uniform topology is that it is not separable anymore. I did not give any proof of this statement, so here it is: all functions of the form $\mathds{1}_{[x,1)}$ are at $\Vert\cdot\Vert_\infty$-distance one of each other. Indeed, if we take $x < y < 1$, then
\[
\left\vert \mathds{1}_{[x,1)}(x) - \mathds{1}_{[y,1)}(x)\right\vert = 1-0 = 1.
\]
This proves that there is an uncountable family of functions which are all at distance one from each other in the uniform topology. Hence, the uniform topology is not separable on $\mathbb{D}_{[0,1]}(\mathbb{R})$.

But let us take another point of view. Perhaps the non separability is not the main problem here. Perhaps it is rather a consequence of an even bigger problem: shouldn't the convergence $\mathds{1}_{\left[\frac{1}{2}-\frac{1}{n},1\right)}\underset{n\to\infty}{\longrightarrow} \mathds{1}_{\left[\frac{1}{2},1\right)}$ hold from an intuitive point of view? That would immediately force these functions to get ``closer" together and prevent the non separability\footnote{at least the one induced by this specific example... but it turns out that this is enough.}.

\begin{figure}[h!tb]
	\centering
	\includegraphics[width=.8\textwidth]{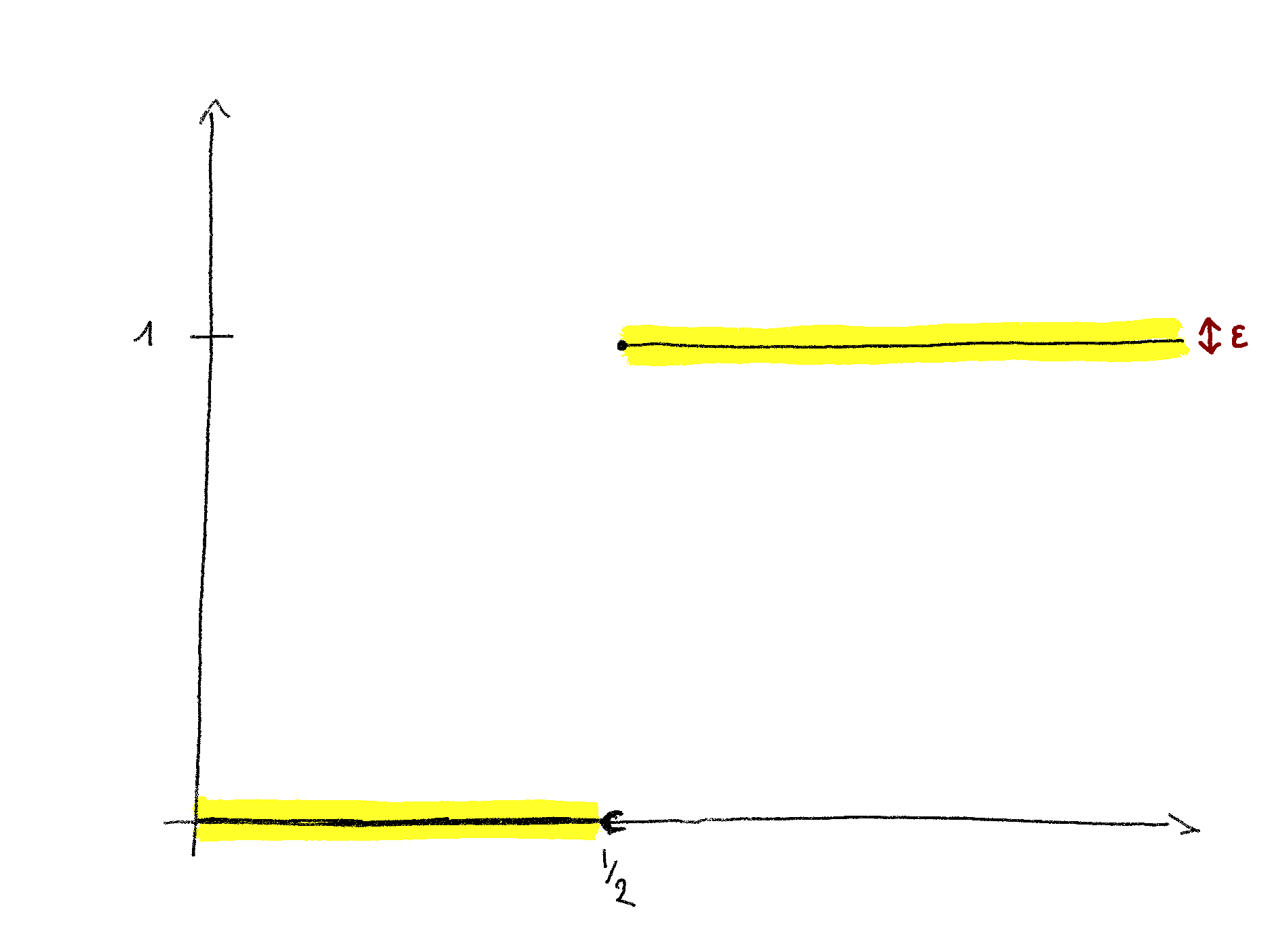}
	\caption{The indicator function $\mathds{1}_{\left[\frac{1}{2},1\right)}$ with its $\epsilon$-tube.}
	\label{fig:eps_tube}
\end{figure}

This insight transforms the problem of finding a topology similar to the uniform topology but preventing non separability into the problem of tweaking the uniform convergence so that this sort of convergence is allowed. To narrow down on what exactly keeps these indicator functions apart, let us have a closer look at what $\epsilon$-balls look like in the uniform topology. Let us take the example of $f = \mathds{1}_{\left[\frac{1}{2},1\right)}$. Then, the $\epsilon$-ball around $f$ contains all functions whose graphs lie in the so-called $\epsilon$-tube around the graph of $f$, see Figure \ref{fig:eps_tube}. This tube forces functions to be ever closer to $f$, but allows them to wiggle a little bit up and down. Note that this corresponds to a \emph{spatial} wiggle. When functions are continuous, that is all perfectly fine, because a wiggle in time can be translated into a wiggle in space. However, when we have a discontinuity, this is not true anymore. As soon as we move the discontinuity a bit to the left or to the right, we necessarily leave the tube and are immediately ``far away" from $f$.

\begin{figure}[h!tb]
	\begin{subfigure}{.49\textwidth}
		\centering
		\includegraphics[width=\linewidth]{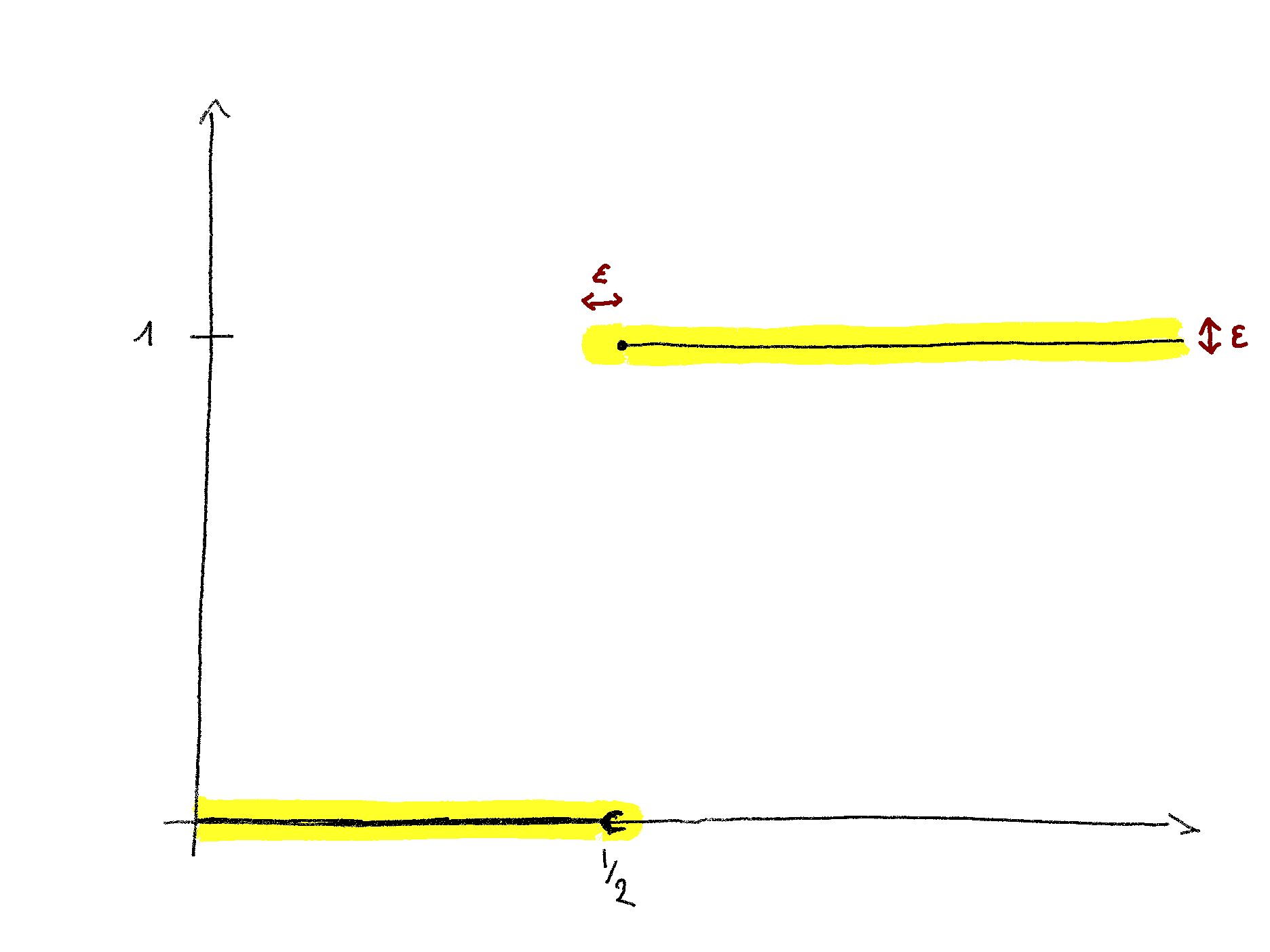}
		\caption{Minimalistic approach}
		\label{subfig:eps_tube_J}
	\end{subfigure}
	\begin{subfigure}{.49\textwidth}
		\centering
		\includegraphics[width=\linewidth]{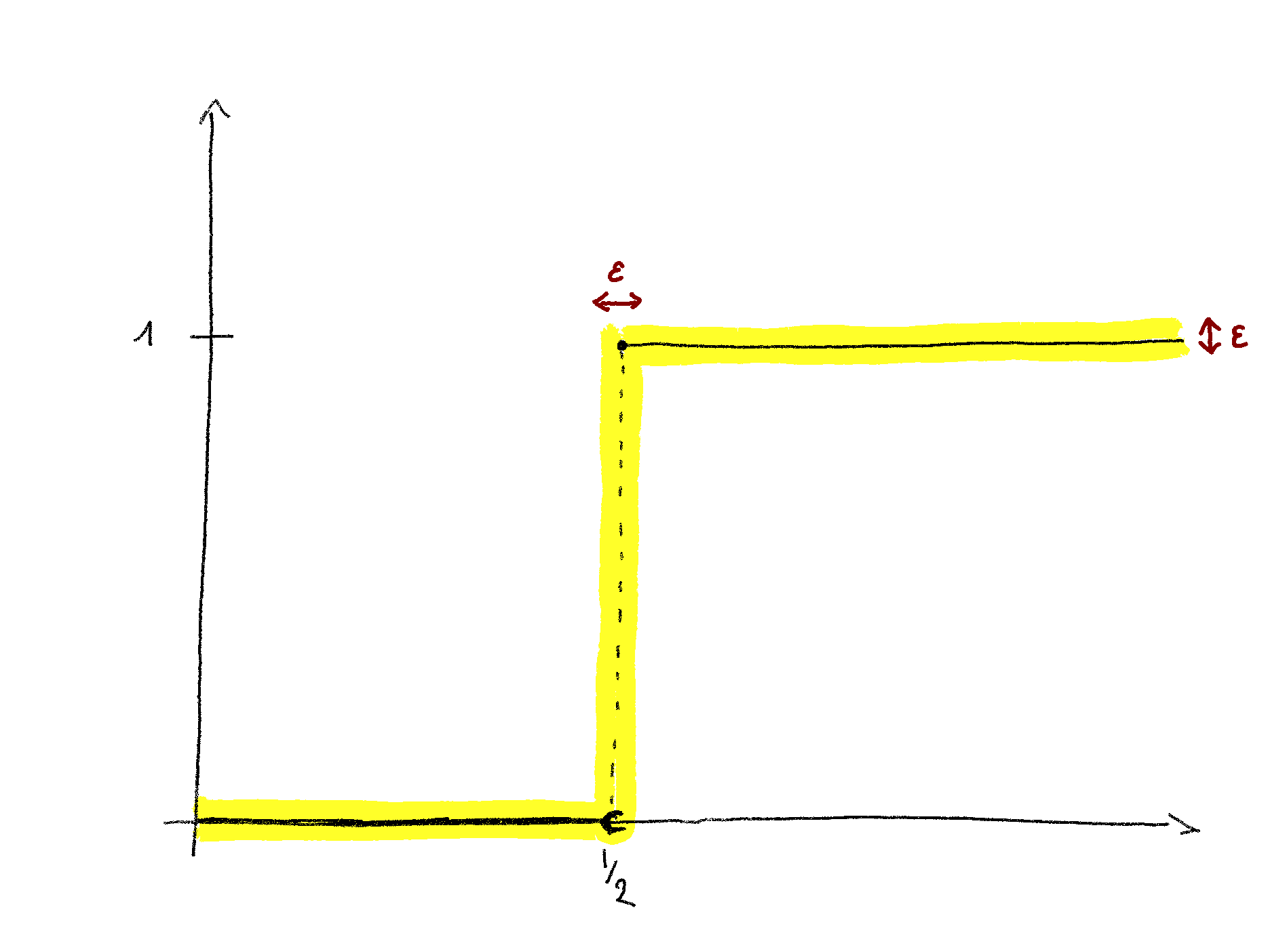}
		\caption{Generous approach}
		\label{subfig:eps_tube_M}
	\end{subfigure}
	\caption{The indicator function $\mathds{1}_{\left[\frac{1}{2},1\right)}$ with modified $\epsilon$-tubes.}
	\label{fig:eps_tube_mod}
\end{figure}

That means that we need to modify the $\epsilon$-tubes to allow for some \emph{temporal} wiggle in addition to the \emph{spatial} wiggle already accounted for. There are two options that may come to mind. The more minimalistic approach would be to extend the $\epsilon$-tube by a little bit at a discontinuity to get $\epsilon$-``gloves", see Figure \ref{subfig:eps_tube_J}. The other, more generous approach, would be to connect the two ends of the $\epsilon$-tube at a discontinuity, see Figure \ref{subfig:eps_tube_M}.

In reality, Skorokhod defined in its foundational paper \cite{Sko} \emph{four} different topologies, a strong and a weak version for each approach. They are now commonly referred to by the rather obscure names $J_1$-, $J_2$-, $M_1$- and $M_2$-topologies. The $J$-topologies arise from the minimalistic approach and the $M$-topologies from the more generous one. 

For all those intimidated by these cryptic names: in the end, only the $J_1$-topology is commonly used and therefore referred to as \emph{the} Skorokhod topology. It appears that the $M$-topologies also have their use in various problems, whereas you will most certainly not encounter the $J_2$-topology at all. That means that this mess reduces to 
\begin{enumerate}[label=\roman*)]
	\item one main ($J_1$-)Skorokhod topology everybody should be familiar with and
	\item a second type of ($M$-)Skorokhod topologies one should have a general idea of.
\end{enumerate}

\begin{figure}[h!tb]
	\centering
	\begin{tikzpicture}
		\node (u) at (0, 2) {$U$};
		\node (j) at (0,0) {$J_1$};
		\node (J) at (-2, -2) {$J_2$};
		\node (m) at (2,-2) {$M_1$};
		\node (M) at (0,-4) {$M_2$};
		
		\draw[->] (u) -- (j);
		\draw[->] (j) -- (J);
		\draw[->] (j) -- (m);
		\draw[->] (m) -- (M);
		\draw[->] (J) -- (M);
	\end{tikzpicture}
	\caption{Relationships between the different Skorokhod topologies. Here $\tau \rightarrow \sigma$ means that $\tau$ is stronger than $\sigma$ in the sense that any sequence converging in $\tau$ does also converge in $\sigma$. $U$ denotes the topology of uniform convergence.}\label{fig:Skorokhod_topologies}
\end{figure}
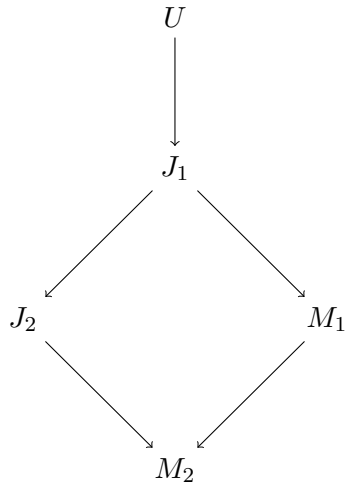

\begin{figure}[h!tb]
	\centering
	\begin{subfigure}{.415\textwidth}
		\centering
		\includegraphics[width=\linewidth]{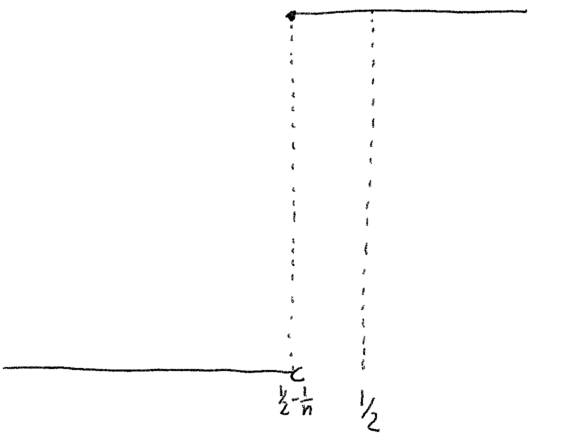}
		\caption{$J_1$}
	\end{subfigure}
	\begin{subfigure}{.4\textwidth}
		\centering
		\includegraphics[width=\linewidth]{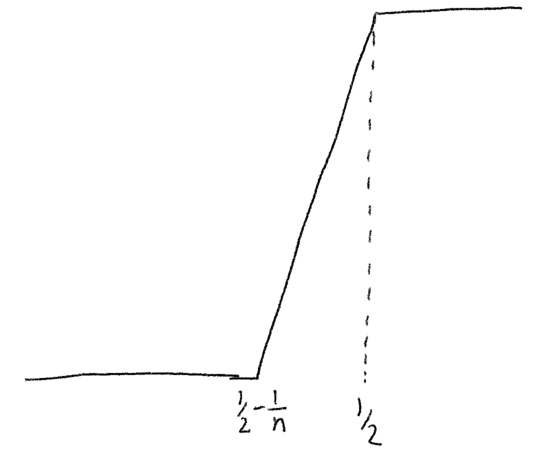}
		\caption{$M_1$ but not $J_2$}
		\label{subfig:M1notJ2}
	\end{subfigure}
	\vspace*{.6cm}
	
	\begin{subfigure}{.407\textwidth}
		\centering
		\includegraphics[width=\textwidth]{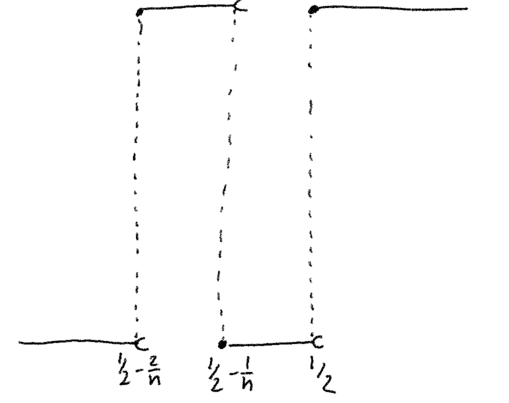}
		\caption{$J_2$ but not $M_1$}
	\end{subfigure}
	\begin{subfigure}{.4\textwidth}
		\centering
		\includegraphics[width=\textwidth]{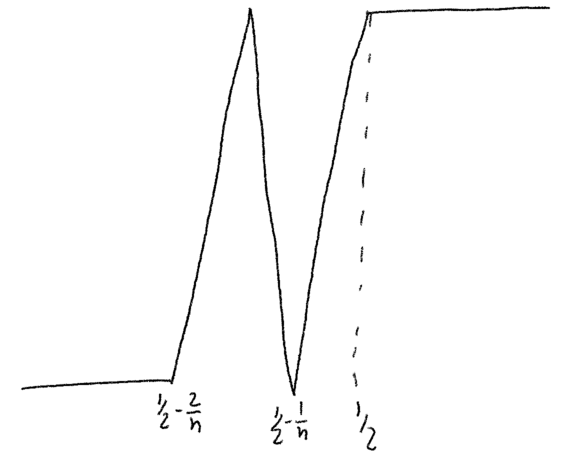}
		\caption{$M_2$}
	\end{subfigure}
	\caption{Examples of sequences converging to $\mathds{1}_{\left[\frac{1}{2},1\right)}$ in the different Skorokhod topologies.}
	\label{fig:diffs_Skorokhod_topologies}
\end{figure}

In this paper, I will only discuss the two main topologies $J_1$ and $M_1$. But before getting into the details, I want to illustrate how these topologies relate to each other. Keeping the above in mind, we should expect the $J$-topologies to be stronger than the $M$-topologies, as the latter allow more functions to be close. However, the situation is a bit more complicated, see Figure \ref{fig:Skorokhod_topologies}. The good news are that \emph{the} Skorokhod topology $J_1$ is stronger than all the other ``new" topologies. That means that whenever a convergence is shown to hold in $J_1$, then it holds in all the other Skorokhod topologies. Figure \ref{fig:diffs_Skorokhod_topologies} illustrates for what extra type of convergence the different topologies allow. These examples are taken from \cite[Figure 11.2]{SPL} and can partially be found already in \cite{Sko}, see also \cite[Limit Theorems for Stochastic Processes]{SW}.

\subsection{The Topologies $J_1$ and $M_1$}\label{ssec:J1andM1}

Now comes the more difficult part of translating our intuition into real definitions. We will start with \emph{the} Skorokhod topology $J_1$ and finish with the $M_1$-topology.

In the minimalistic setting, we want to allow for some temporal wiggle, without connecting the two ends of the $\epsilon$-tube across a discontinuity. To reformulate this mathematically, we will perform a time change. By a \emph{change of time} I mean that we take a strictly increasing bijection $\lambda :[0,1]\rightarrow [0,1]$ and consider $f\circ \lambda$ instead of $f$. Naturally, we are only interested in parametrisations that are ``close" to the unitary flow of time, corresponding to the trivial parametrisation $id: t\mapsto t$. In other words, we need to penalise parametrisations which are ``too far away" from $id$. This leads to the following definition of distance:
\[
d_{J_1}(f,g) := \inf_{\lambda} \Big{\{} \Vert f\circ \lambda - g\Vert_\infty  + \Vert \lambda - id\Vert_\infty\Big{\}},
\]
where the infimum is taken over all increasing bijections on $[0,1]$. It can be shown that $d_{J_1}$ is a metric and it is usually referred to as \emph{Skorokhod metric}, see \emph{e.g.}~\cite[Section 12]{Billingsley} and one defines $J_1$ to be the topology induced by this metric.

\begin{figure}[h!tb]
	\centering
	\includegraphics[width=.8\textwidth]{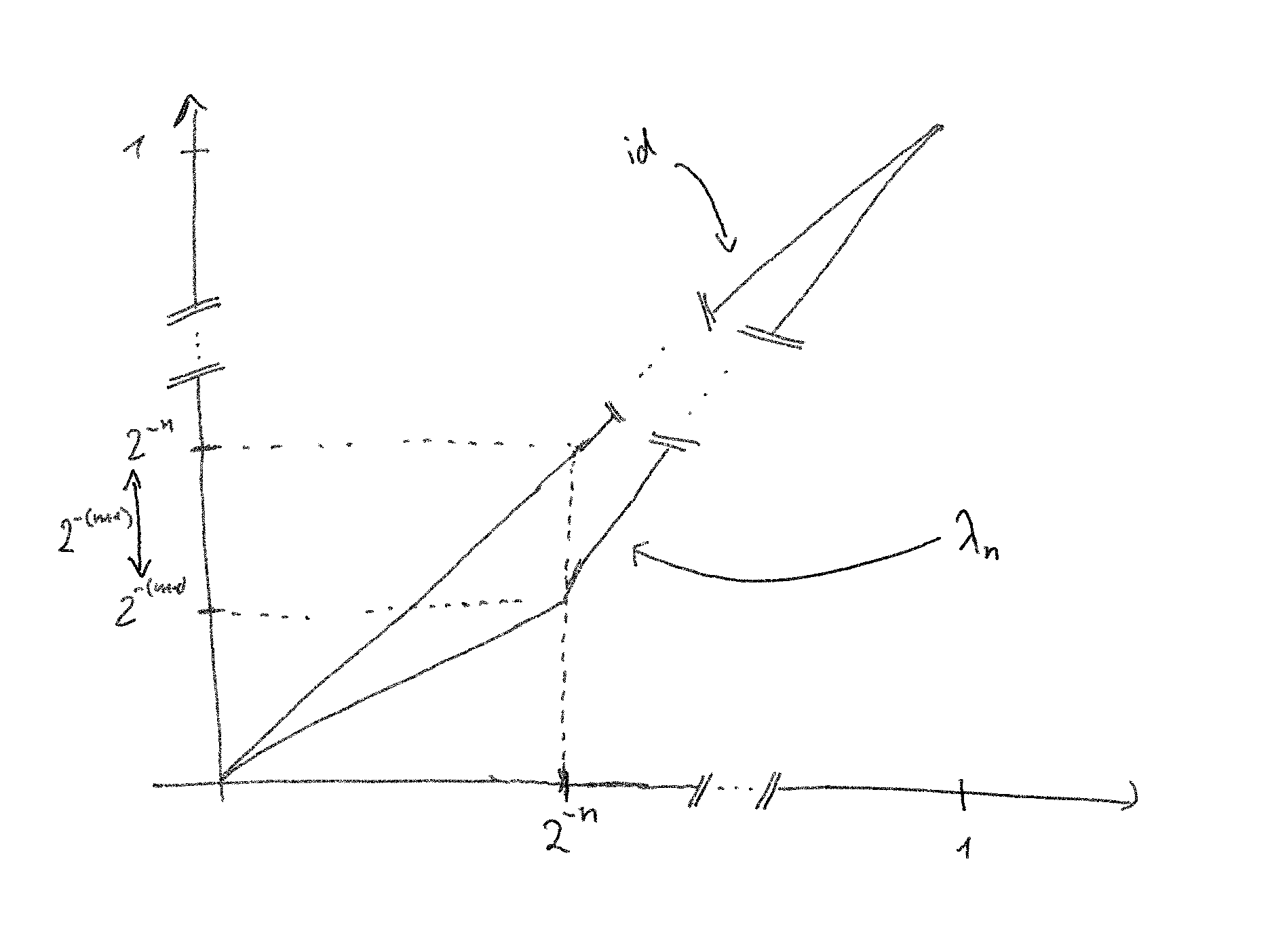}
	\caption{The change of time $\lambda_n$ compared to the usual time flow $id$.}\label{fig:metric_incomplete}
\end{figure}

Unfortunately, it turns out that this is actually a bad metric in the sense that it is not complete. Consider the following example from \cite[Example 12.2]{Billingsley}. Take the indicator functions $f_n := \mathds{1}_{[0, 2^{-n})}$ and define the change of time $\lambda_n$ as the linear interpolation of the three points $(0,0)$, $(2^{-n}, 2^{-(n+1)})$ and $(1,1)$, see Figure \ref{fig:metric_incomplete}. One easily checks that the change of time is such that $f_n = f_{n+1}\circ\lambda_n$ leading to $\Vert f_{n+1}\circ\lambda - f_n\Vert_\infty = 0$. To bound the penalty $\Vert \lambda_n - id\Vert_\infty$ on the change of time, note that they differ maximally at $x = 2^{-n}$ yielding $\Vert \lambda_n - id\Vert_\infty= \vert 2^{-n} - 2^{-(n+1)}\vert = 2^{-(n+1)}$. This gives
\[
d_{J_1}(f_{n+1}, f_n) \leq \Vert f_{n+1}\circ\lambda_n - f_n\Vert_\infty + \Vert \lambda_n - id\Vert_\infty = 2^{-(n+1)}.
\]
In particular, this error is summable and we conclude that $(f_n)_n$ is a Cauchy sequence. Since $f_n(x)$ converges to $0$ for all $x\in(0,1)$, the only possible limit is the null function $f=0$. However, whatever change of time we apply, the null function doesn't change. Hence,
\[
d_{J_1}(f, f_n) = \Vert f_n\Vert_\infty = 1.
\]
That means that although $(f_n)_n$ is Cauchy, it does not converge.

The problem of the Skorokhod metric is quite subtle and lies within the penalty we put onto the parametrisation $\lambda$: we measure the absolute distance between $\lambda$ and $id$. However, it is better to think of parametrisations as a modified \emph{flow of time}. In this sense, it would be better to measure the difference in flow \emph{speed}. In other words, we want parametrisations with a nearly constant speed
\[
\dfrac{\lambda(t) - \lambda(s)}{t-s}\approx 1.
\]
In the above example, the slope of $\lambda_n$ will never converge to $1$ (not even pointwise), as
\[
\dfrac{\lambda_n(2^{-n}) - \lambda(0)}{2^{-n}} = 2^n\cdot 2^{-(n+1)} = \dfrac{1}{2}
\]
for every $n\geq 1$. To fix this, we introduce the new penalty
\[
\Vert \lambda\Vert_\circ := \sup_{t\neq s} \left\vert \log\left(\dfrac{\lambda(t) - \lambda(s)}{t-s}\right)\right\vert,
\]
leading to the modified metric
\[
\hat{d}_{J_1}(f,g) := \inf_\lambda \Big{\{} \Vert f \circ \lambda - g\Vert_\infty + \Vert \lambda\Vert_\circ\Big{\}}.
\]

It turns out that both metrics are equivalent, \emph{i.e.}~induce the same topology. However, this modified metric is complete! For this reason, this metric is sometimes called \emph{Skorokhod metric} instead of the previous one, leaving the original metric without any special name. It can be shown that $J_1$ is separable (see \emph{e.g.}~\cite[Theorem 12.2]{Billingsley}), meaning that $J_1$ satisfies the \eqref{eq:polish}. 

Recalling that the uniform topology on continuous functions does not care about small temporal distortions, one easily verifies that $J_1$ also satisfies the \eqref{eq:extension}, \emph{i.e.} if a sequence of continuous functions converges in $J_1$, then the convergence is uniform (and conversely).

The only thing we still need is a good description of compact sets, \emph{i.e.}~that $J_1$ satisfies the \eqref{eq:aa}. Fortunately, there is indeed a result similar to the Arzelà--Ascoli Theorem! The only thing we need to adapt is the definition of the modulus of continuity so that it ignores jump discontinuities. This is achieved by allowing the function to jump at a finite number of points:
\[
\omega'_\delta(f) := \inf_{\{t_i\}_{i=0}^{v}} \;\,\sup_{1\leq i\leq v} \;\,\sup_{s,t\in [t_{i-1},t_i)} \vert f(t) - f(s)\vert,
\]
where the infimum is taken over all finite partitions $0=t_0 < \dots < t_v = 1$ of $[0,1]$ of finite size $v\geq 1$. Note that the innermost supremum only ranges over the right open interval $[t_{i-1}, t_i)$, allowing for jumps at times $t_i$. To distinguish it from the ``real" modulus of continuity, I refer to it as \emph{modulus of continuity type function}.
\begin{theo}[Compactness in $J_1$, see \emph{e.g.}~{\cite[Theorem 12.3]{Billingsley}}]\label{theo:compactness_J1}
A set\\ $K\subseteq \mathbb{D}_{[0,1]}(\mathbb{R})$ is relatively compact in $J_1$ if and only if
\begin{enumerate}[label=\roman*)]
	\item the set is uniformly bounded in the sense that
	\[
	\sup_{f\in K} \Vert f\Vert_\infty < +\infty;
	\]
	\item the modulus of continuity type function vanishes uniformly over $K$:
	\[
	\lim_{\delta\downarrow 0} \sup_{f\in K} \omega_\delta'(f) = 0.
	\]
\end{enumerate}
\end{theo}

Note that it is not enough to have a uniform bound of $\vert f(0)\vert$ as before! This was possible in the case of continuous functions, because we also imposed uniform equicontinuity. Since we now allow for jumps, we need to strengthen this condition to a uniform bound on the entire interval.

That means that $J_1$ satisfies all of the constraints that are really important to us. The only thing we are left to check is the \eqref{eq:BONUS} of strengthened convergence whenever the limiting function is continuous. From the heuristic, it is conceivable that this holds for $J_1$ as we only allow for small distortions in time. This idea can be made rigorous by noting that we may shift the change of time onto the limit: take $f_n \to f$ in $J_1$ with $f$ continuous. Then there is a sequence $(\lambda_n)_n$ of time changes such that
\[
\lim_n \Big{(} \Vert f\circ\lambda_n - f_n\Vert_\infty + \Vert \lambda_n\Vert_\circ\Big{)} = 0.
\]
Since $f$ is continuous and $\Vert \lambda\Vert_\circ\to 0$, one has $\Vert f\circ\lambda_n - f\Vert_\infty \to 0$. This in turn implies that $\Vert f - f_n\Vert_\infty \to 0$.
A extension of this statement can be found \emph{e.g.}~in \cite[Proposition 3.6.5]{EK86} and the preceding comment.

We finish with the $J_1$ topology by pointing out that there is an immediate drawback in satisfying the \eqref{eq:extension}: continuous functions cannot converge to discontinuous functions as $\mathcal{C}([0,1])$ is closed in the uniform topology! This explains why we cannot expect to see convergences of the type shown in Figure \ref{subfig:M1notJ2}. \\

Let us now come to the $M_1$-topology. Recall that the idea was to generously extend the $\epsilon$-tube across the discontinuity, see Figure \ref{subfig:eps_tube_M}. To make this more rigorous, we will work with the so-called \emph{completed graph}
\[
\Gamma(f) := \{ (t,z)\in [0,1]\times \mathbb{R}\;:\; z = \alpha f(t-) + (1-\alpha)f(t)\text{ for some }\alpha\in[0,1]\}
\]
of $f\in \mathbb{D}_{[0,1]}(\mathbb{R})$ containing the graph of $f$ together with the straight lines connecting the two ends of a jump discontinuity. Here, I use the notation $f(t-)$ to denote the left limit of $f$ in $t$ which exists by definition of $\mathbb{D}_{[0,1]}(\mathbb{R})$. For example, in Figure \ref{subfig:eps_tube_M}, the completed graph corresponds to the graph together with the dotted line.

From here on, the idea is very similar to the one used for the $J_1$-topology. To allow for some freedom, we use again parametrisations. The only difference is that we will now parametrise the completed graph, \emph{i.e.}~we will take non decreasing functions $(\lambda, \rho):[0,1)\rightarrow \Gamma(f)$ which are onto. Here, $\lambda$ is the temporal component and $\rho$ is the spatial component. Note also that we use the intuitive order on $\Gamma(f)$ to define monotonicity:
\[
(t_1,z_1)\leq (t_2,z_2)\quad\text{iff}\quad \Big{(} t_1 < t_2\quad\text{ or }\quad t_1 = t_2\text{ and }\vert z_1 - f(t-)\vert \leq\vert z_2 - f(t-)\vert \Big{)}.
\]
In words, we use the temporal ordering coming from drawing the completed graph from left to right without lifting the pencil. We then define the distance between two functions $f$ and $g$ through
\[
d_{M_1}(f,g) := \inf_{(\lambda_f,\rho_f), (\lambda_g,\rho_g)} \Big{\{} \Vert \lambda_f - \lambda_g\Vert_\infty + \Vert \rho_f - \rho_g\Vert_\infty\Big{\}}
\]
as the minimal distance between any two parametric representations of $f$ and $g$. Again, one checks that $d_{M_1}$ indeed is a metric and one defines $M_1$ to be the induced topology.

The advantage of this metric is that we get the additional convergence illustrated in Figure \ref{subfig:M1notJ2}. Symmetrically, the drawback is that $M_1$ does not verify the \eqref{eq:extension}. Except from this, all other constraints are satisfied: $M_1$ defines a Polish topology\footnote{Again, the metric we defined is \emph{not} complete. An equivalent complete metric is defined in \cite[Section 12.8]{SPL}.} on $\mathbb{D}_{[0,1]}(\mathbb{R})$, see \emph{e.g.}~\cite[Theorem 12.8.1]{SPL} and also satisfies the \eqref{eq:BONUS} property of strengthened convergence whenever the limit is continuous. The description of compact sets is a bit more difficult, but it is still possible. Before we define the new modulus of continuity we will need here, one should keep in mind that $M_1$ is weaker than $J_1$. That means in particular that the above criterion for compactness still is \emph{sufficient}. Then again, in most of the cases if the above applies, one usually directly works with the $J_1$-topology.

Define the new modulus of continuity type function
\[
\omega_\delta''(f) := \sup_{t\in [0,1]} w(f,t,\delta)
\]
with
\[
w(f,t,\delta) := \sup_{0\vee t-\delta \leq t_1 < t_2 < t_3< t+\delta\wedge 1}\;\, \inf_{z\in [f(t_1),f(t_3)]} \vert f(t_2) - z\vert.
\]
Using the notation $d(x,A)$ for the distance between $x$ and the set $A$, this can be written more compactly as
\[
w(f,t,\delta) := \sup_{0\vee t-\delta \leq t_1 < t_2 < t_3< t+\delta\wedge 1} d\big(f(t_2), [f(t_1),f(t_3)]\big).
\]
A small modulus of continuity type function ensures that on small intervals, the graph is close to straight lines. This is a way to exclude oscillations, but allows for ever steeper slopes, see Figure \ref{subfig:M1notJ2}.

\begin{theo}[Compactness in $M_1$, see \emph{e.g.}~{\cite[Theorem 12.12.2]{SPL}}]
A set $K\subseteq \mathbb{D}_{[0,1]}(\mathbb{R})$ is compact in $M_1$ if and only if
\begin{enumerate}[label=\roman*)]
	\item the set is uniformly bounded in the sense that
	\[
	\sup_{f\in K} \Vert f\Vert_\infty < +\infty;
	\]
	\item the oscillations vanish uniformly over $K$:
	\[
	\left\{\begin{array}{rll}
	\lim_{\delta\downarrow 0} \sup_{f\in K} &\omega_\delta''(f) &= 0\\
	\lim_{\delta\downarrow 0} \sup_{f\in K} &\vert f(\delta) - f(0)\vert &= 0\\
	\lim_{\delta\downarrow 0} \sup_{f\in K}&\vert f(1-) - f(1-\delta)\vert &= 0
	\end{array}\right..
	\]
\end{enumerate}
\end{theo}

The above description seems to differ from the short characterisation in Theorem \ref{theo:compactness_J1}, but a similar description for compact sets w.r.t.~$J_1$ can be found \emph{e.g.}~in \cite[Theorem 12.4]{Billingsley}.\\

We finish this section with a small overview of which topologies have which properties. For completeness, I also include the topologies $J_2$ and $M_2$. These additional properties follow columnwise from \cite[Theorem 11.6.6]{SPL}; \cite[Limit Theorems for Stochastic Processes, Section 2.7]{SW} and \cite[Theorems 12.12.2]{SPL}; a similar argument for $J_2$ as in Section 3 and for $M_2$ the fact that it is weaker than $M_1$; and \cite[Corollary 12.11.1]{SPL} together with the fact that $J_2$ is stronger than $M_2$.
\begin{center}
	\begin{tabular}{| c || c | c | c | c | c |}\hline
	& \small\hyperref[eq:polish]{Polish} & \small\hyperref[eq:prokhorov]{SepPrkhorov} & \small\hyperref[eq:aa]{Arzelà--Ascoli} & \small\hyperref[eq:extension]{Extension Property} & \small\hyperref[eq:BONUS]{Bonus Property} \\\hline\hline
	$J_1$ & \color{darkgreen}Yes & \color{darkgreen}Yes & \color{darkgreen}Yes & \color{darkgreen}Yes & \color{darkgreen}Yes \\\hline
	$M_1$ & \color{darkgreen}Yes & \color{darkgreen}Yes & \color{darkgreen}Yes & \color{red}No & \color{darkgreen}Yes \\\hline
	$J_2$ & \color{red}? & \color{darkgreen}Yes & \color{darkgreen}Weakly Yes & \color{darkgreen}Yes & \color{darkgreen}Yes\\\hline
	$M_2$ & \color{red}? & \color{darkgreen}Yes & \color{darkgreen}Weakly Yes & \color{red}No & \color{darkgreen}Yes\\\hline
	\end{tabular}
\end{center}

\subsection{Convergence of Stochastic Processes in $J_1$ and $M_1$}

In this section, I simply restate the previous compactness results in terms of tightness of stochastic processes on $\mathbb{D}_{[0,1]}(\mathbb{R})$. I will furthermore state the generalisation of Corollary \ref{corol:convergence_law_equiv} to the topologies $J_1$ and $M_1$. This is not immediate, because the projections $\pi_t: f\mapsto f(t)$ are not continuous anymore: already in $J_1$, the projection $\pi_t$, $t\in (0,1)$, is continuous in $f$ if and only if $f$ is continuous in $t$. (The projections $\pi_0$ and $\pi_1$ are always continuous.) That means that we cannot hope to have the convergence of all finite-dimensional distributions.

Instead, we define the set of all continuity points of a stochastic process $X$ by
\begin{align*}
T_X &:= \{ t\in (0,1)\;:\; \mathbb{P}(X \text{ is continuous at $t$}) = 1\}\cup \{0,1\} \\
&\;= \{t\in [0,1]\;:\; \mathbb{P}(\pi_t\text{ is continuous in $X$}) = 1\}.
\end{align*}
Note that $0$ is always a continuity point of any stochastic process $X$, but $X$ may be discontinuous in $t=1$. As such, it would be better to say that $T_X$ is the set of times $t\in [0,1]$ such that the projection $\pi_t: f\mapsto f(t)$ is a.s.~continuous.
It turns out to be enough to check convergence on $T_X$, which is almost surely dense in $[0,1]$, see \cite[Section 13]{Billingsley}.

\begin{theo}[see \emph{e.g.}~{\cite[Theorem 11.6.6]{SPL}}]\label{theo:tight+finite_dim_Skorokhod}
A sequence of \emph{càdlàg} processes $(X_n)_{n\geq 1}$ converges in law to a \emph{càdlàg} process $X$ w.r.t.~either $J_1$ or $M_1$ if and only if $(X_n)_{n\geq 1}$ is tight in the respective topology and all finite dimensional distributions at times $t_i\in T_X$ converge to those of $X$.
\end{theo}

Let us now turn to tightness criteria.

\begin{theo}[Tightness in $J_1$, see \emph{e.g.}~{\cite[Theorem 13.2]{Billingsley}}]
A sequence of \emph{càdlàg} processes $(X_n)_{n\geq 1}$ is tight w.r.t.~$J_1$ if and only if
\begin{enumerate}[label=\roman*)]
	\item it holds that
	\[
	\lim_{C\to+\infty} \limsup_n \mathbb{P}\left[ \Vert X_n\Vert_\infty \geq C\right] = 0;
	\]
	\item for every $\epsilon > 0$, it holds that
	\[
	\lim_{\delta\downarrow 0} \limsup_n \mathbb{P}\left[ \omega_{\delta}'(X_n) \geq \epsilon \right] = 0.
	\]
\end{enumerate}
\end{theo}

\begin{theo}[Tightness in $M_1$, see \emph{e.g.}~{\cite[Theorem 12.12.3]{SPL}}]
A sequence of \emph{càdlàg} processes $(X_n)_{n\geq 1}$ is tight w.r.t.~$M_1$ if and only if
\begin{enumerate}[label=\roman*)]
	\item it holds that
	\[
	\lim_{C\to+\infty} \limsup_n \mathbb{P}\left[ \Vert X_n\Vert_\infty \geq C\right] = 0;
	\]
	\item for every $\epsilon > 0$, it holds that
	\[
	\begin{cases}
	\lim_{\delta\downarrow 0} \limsup_n \mathbb{P}\left[ \omega_{\delta}''(X_n) \geq \epsilon \right] &= 0\\
	\lim_{\delta\downarrow 0} \limsup_n \mathbb{P}\left[ \vert X_n(\delta) - X_n(0)\vert > \epsilon\right] &= 0\\
	\lim_{\delta\downarrow 0} \limsup_n \mathbb{P}\left[ \vert X_n(1-) - X_n(1-\delta)\vert > \epsilon\right] &= 0\\
	\end{cases}.
	\]
\end{enumerate}
\end{theo}

For sufficient criteria of tightness in $J_1$ in a very general setting, one may refer to \cite[Chapter 3]{EK86}, or \cite[Chapter 1]{Etheridge} for a less general approach. Convergence criteria for $M_1$ and $M_2$ may be found in \cite[Chapter 12]{SPL}.

\section{Extending $J_1$ and $M_1$ to More General Skorokhod Spaces}

There are some generalisations we are interested in. The first concerns the time interval: obviously, there is no problem with substituting $[0,1]$ with some other finite time interval $[0,T]$, but what about the entire half line $[0,+\infty)$? Secondly, we would like to be able to consider the space $\mathbb{D}_{[0,T]}(\mathbb{R}^k)$ and compare its topology with the product topology $\mathbb{D}_{[0,T]}(\mathbb{R})^k$. Finally, we would like to generalise the topology to a more general range space $E$ to get Skorokhod topologies on $\mathbb{D}_{[0,T]}(E)$.

This section is relatively short as I will simply point out possible restrictions and the necessary modifications needed to generalise the topologies.

\subsection{Extending Time}

As mentioned above, there is no difficulty in extending the topology to $\mathbb{D}_{[0,T]}(\mathbb{R})$ for any finite time horizon $T > 0$. One simply replaces the $1$s in all definitions by $T$s.

Once we can define the topology on any finite time horizon, we would like to extend it to $\mathbb{D}_{[0,+\infty)}(\mathbb{R})$ in the usual way: a sequence $(f_n)\subseteq \mathbb{D}_{[0,+\infty)}(\mathbb{R})$ converges to $f\in \mathbb{D}_{[0,+\infty)}(\mathbb{R})$ if and only if all restrictions to time intervals of the form $[0,T]$ converge to the restrictions of $f$ to these intervals. 
However, this approach is doomed: the sequence given by $f_n := \mathds{1}_{[1+1/n, +\infty)}$ does not converge to $\mathds{1}_{[1,+\infty)}$ in $\mathbb{D}_{[0,1]}(\mathbb{R})$! The problem is that the endpoint of the time interval is different from all other points. Since this problem disappears whenever it is a continuity point of the limit function, one remedy is to say that $f_n$ converges to $f$ if and only if all restriction to the time intervals of the form $[0,T]$, with $T$ a continuity point of $f$, converge to the restrictions of $f$ to these intervals. A different approach is demonstrated in \cite[Section 16]{Billingsley}.

It is possible to define a metric which induces this topology, see \emph{e.g.}~\cite[Sections 3.3 and 12.9]{SPL} or \cite[Section 16]{Billingsley}. One checks that $\mathbb{D}_{[0,+\infty)}(\mathbb{R})$ is again a Polish space, see \emph{e.g.}~\cite[Theorem 16.3]{Billingsley} in the case of $J_1$. By our definition, compactness can be checked by restricting to compact time intervals and using the compactness criteria discussed above.

\subsection{Product Skorokhod Topologies}

The next step is to generalise the topologies to processes with values in $\mathbb{R}^k$. The easiest solution is to wait for the next section and take the range space $\mathbb{R}^k$. However, there is a second way we can get a topology on $\mathbb{D}_{[0,1]}(\mathbb{R}^k)$ by using the fact that
\[
\mathbb{D}_{[0,1]}(\mathbb{R}^k) \equiv \left(\mathbb{D}_{[0,1]}(\mathbb{R})\right)^k
\]
in the sense of a bijection. From the point of view of topologies, this implies that we may endow $\mathbb{D}_{[0,1]}(\mathbb{R}^k)$ with the product topology on $\left(\mathbb{D}_{[0,1]}(\mathbb{R})\right)^k$. It turns out that this topology differs from the topology we would get from the next section. More precisely, the product topology is \emph{weaker}. For this reason, we speak of the strong topologies obtained by viewing $\mathbb{R}^k$ as the range space and the weak topologies obtained as product topologies. Intuitively, the product topology allows for different time parametrisations in every component whereas in the strong topology, the same time change is used for all components.

A detailed discussion of the weak $M$-topologies can be found in \cite[Chapter 12]{SPL}. Unfortunately, I have not found much literature on the weak $J$-topologies. It might be that they are less commonly used. One use may be found \emph{e.g.}~in \cite[Proof of Lemma 4.3]{EK19}, see also \cite[Theorem 11.5.1]{SPL}.

\subsection{Towards General Range Spaces}

Let $E$ be a general Polish space. When going back to the definition of $J_1$ in Section \ref{ssec:J1andM1}, we only need to replace $\Vert f\circ\lambda - g\Vert_\infty$ by
\[
\sup_{t\in [0,1]} d_E\Big{(}(f\circ\lambda)(t), g(t)\Big{)},
\]
where $d_E$ is a complete metric on $E$. It turns out that everything else goes through without any problem. To generalise the compactness result, we only need to adjust again the definition of the modulus of continuity to
\[
\omega_\delta'(f) := \inf_{\{t_i\}_{i=0}^v} \;\,\sup_{1\leq i\leq v} \;\,\sup_{s,t\in [t_{i-1},t_i)} d_E\Big{(} f(t), f(s)\Big{)}.
\]
One can verify that $J_1$ preserves all the good properties as long as $E$ is Polish, see \emph{e.g.}~\cite[Chapter 3]{EK86}.\\

One would hope that this is also possible with the $M_1$-topology. However, its definition relies on the fact that we can define the straight line between two points in $E$. In other words, we need an additive structure to define $M_1$, \emph{i.e.}~we can generalise $M_1$ only to Banach spaces. In this setting, the interval $[f(t_1), f(t_2)]$ has to be interpreted as the line 
\[
\{ (1-\alpha)f(t_1) + \alpha f(t_2)\;:\; \alpha\in[0,1]\}
\]
from $f(t_1)$ to $f(t_2)$. This restriction is another reason why $J_1$ has become the more prevalent topology.

\printbibliography

\end{document}